\documentclass[12pt]{article}

\usepackage{amsmath}
\usepackage{amsfonts, amssymb}
\usepackage[cp1251]{inputenc}
\usepackage[russian]{babel}

\def\P{{\bf {P}}}

\def\E{{\bf {E}}}

\def\a{{\alpha}}
\def\d{{\delta}}
\def\si {\sigma}

\def\la {{\lambda}}

\def\ep {{\varepsilon}}
\def\ga {{\gamma}}
\def\Ga {{\Gamma}}
\def\Dl {{\Delta}}

\def\Ka {{\varkappa}}

\begin{document}

\centerline {\bf Rozovsky L.}

\medskip

{\bf Small ball probabilities for certain gaussian fields.}

\bigskip

{\bf 1. Introduction and results.}

\bigskip 

We study the behavior of the tail probabilities $\P(V^2<r)$ as
$r\to 0$, where $V^2$ is defined by the following double sum
$$
V^2 =\pi^{-4}\,\sum\limits_{i,j\ge 1} \Big((i+b)\,(j +
\d)\Big)^{-2}\,\xi_{ij}^2,\eqno (1.1)
$$
where $ \{\xi_{ij}\} $ are independent standard normal random
variables, and $b$ and $\d$ are constants: $b>-1,\ \d>-1$.

According to the well-known Karhunen –- Lo\`eve expansion,  the
sum in (1.1)  can be considered as the squared $ L_2 $-norm of a
2--parameter Gaussian random field, which has the covariance of a
``tensor product'' type (see,  for instance, [6] and [7]). The
classical examples of such fields are the standard Brownian sheet
$W(t, s)$, i.e. $V^2=\int_0^1 \int_0^1 W^2(t, s)\,dt\,ds$, the
Brownian pillow and the pillow-sleep, or the Kiefer field ($b =
\d=-1/2$,\ $b = \d=0$ and $b = 0$, $\d=-1/2$, respectively).

We remark that small deviations for Gaussian random fields in the
Hilbert norm have not been studied as extensively as their
one-parameter counterparts. In [4] -- [9], the logarithmic $ L_2 $
small ball asymptotics for some two-parameter and multi-parameter
Gaussian random fields were obtained. Concerning sharper
estimates, the only work we know is [5], where small deviations in
$L_2$--norm for the Brownian sheet and its integrated counterparts
were thoroughly examined. Therefore, a more detailed study of the
tail probabilities $\P(V^2<r)$ is of interest, even though the sum
(1.1) has a special form.

\medskip

In order to formulate our results, we need to introduce some
notation.

\medskip

Let us introduce functions $B_i = B_i(b, \d),\
b>-1,\,\d>-1,\,(1\le i\le 3)$ as follows:
$$\gathered
B_1 = 1 + \Psi(1+b) + \Psi(1+\d),\quad B_2
= (1 + 2\,b)\,(1 + 2\,\d)/2,\\
B_3 = (1 + \d + b)\,\log{(2\,\pi)} - (1 + 2\,\d)\,\ln{\Ga(1 + b)}
 -(1 + 2\,b)\,\ln{\Ga(1 + \d)}.\endgathered\eqno (1.2)
$$
Here $\Psi (x) $ is the logarithmic derivative of the gamma --
function $ \Ga (x) $. Let us set
$$
B = \pi^{-1}\, (B_1 + 2\, \ln\pi),\
 C =B_2,\  D = B_3 - 2\,B_2\, \ln\pi. \eqno (1.3)
$$

\medskip

Next, we introduce polynomials $\pi_k = \pi_k(s),\ k\ge 0,$ such
that $\pi_0 = 1,\ \pi_1 = s$, and for each $ k \ge 2 $ the
function $\pi_{k+1}$ is given by the following recursive relation
$$
\pi_{k+1} = \pi_k - \frac{k-1}{2 k}\,  \sum_{l= 1}^k
\pi_l\,\pi_{k+1-l} ,\ k\ge 1. \eqno (1.4)
$$
So, $ \pi_2 = s,\, \pi_3 = s - s^2/2, \, \pi_4 = s - 3 s^2/2 +
s^3/3$.

\medskip

Remark 1. The polynomials $ \pi_m (s) $ can be determined
explicitly:
$$
\pi_m(s) = \sum_{k=1}^{m-1} s^k\, a_{k-1, m-k},\quad m\ge 2.\eqno
(1.5)
$$
Here the coefficients $a_{jl} = s(j+l, l)/(j+1)!$, and the values
$ s (n, k) $ are the Stirling numbers of the first kind, i.e. (see
[2, \, V]):
$$
\frac{1}{k!}\,\ln^k(1+t) =  \sum_{n\ge k} s(n,k)\,\frac{t^n}{n!}.
$$
In particular, the coefficients in (1.5) at $ s $ and $s^{m-1}$
are $1$ and $ (- 1) ^ m / (m-1) $, respectively.

Moreover, the coefficients $ a_ {jl} $ may be calculated using the
following formulas: $a_{0l}=1$, $a_{j0}\Big|_{j\ge 1} = 0$\quad
and
$$
a_{jl} = \frac{(j+l)!}{(j+1)!} \,(-1)^j\,\sum_{k=1}^{\min{(j, l)}}
\sum \frac{1}{(l-k)!} \prod\limits_{m=1}^j
\frac{(m+1)^{-k_m}}{k_m!},\quad j, l\ge 1,\eqno (1.6)
$$
where the summation is taken over all non-negative integer
solutions $(k_1,\dots,k_j)$ of  the equations  $1\cdot k_1 +\dots
+j \cdot k_j = j$ and $k_1 +\dots + k_j=k$, and as a result,
$$
|a_{jl}|\le  \frac{1}{j+1} \,4^{j+l-1}.\eqno (1.7)
$$

\medskip

We check Remark 1 in Section 2 (see Lemma 2 and below).

\medskip

For each $ r \in (0, 1), $ we introduce the following notation:
$$
\bar s = \ln{|\ln{r}|}  -  \Psi(1+b) - \Psi(1+\d) -
\ln{(2\pi^3)},\quad \bar \pi_m = \pi_m(\bar s),\  m\ge 0. \eqno
(1.8)
$$
Here $\Psi$ and $\pi_m$ occur respectively in (1.2) and (1.4).
Further, let us define the function $ \ep (r) $ as follows:
$$
\ep(r) = \sum_{m\ge 1} p_m(\bar s)\,|\ln{r}|^{-m},\eqno (1.9)
$$
where
$$
p_1(\bar s) = 2\,\bar s - 2,\quad p_2(\bar s) = \bar s^2,\quad
p_m(\bar s) = -\frac{2}{m-2}\,(\bar \pi_m -\bar\pi_{m-1}),\ m\ge
3.
$$


\medskip

Now, using the notation in (1.2), (1.3), (1.8) and (1.9), we are
able to formulate the main result.

\medskip

{\bf Theorem 1.} {\it Let $\Ka(r) = (8 \pi^2 r)^{-1}\, \ln^2 r,$ $
c_0 = \Big(\pi\,2^{ C/2}\,e^{ D}\Big)^{-1/2}$,\ $c_1 = C/4,\ c_2 =
-( C + 2)/4$, where $C$ and $D$ are introduced in $(1.3)$.  Then
as $r\to 0$
$$
\P(V^2 < r) = c_0\,r^{c_1}\,\Ka^{c_2}(r)\, e ^ {- \Ka(r)\,(1
+\ep(r))} \Big (1 +  O\,\Big(\frac{\ln{|\ln r|}}{|\ln r|}\Big)
\Big), \eqno (1.10)
$$
and, in addition, 
$$
\frac{d}{dr} \P(V^2 < r) = \frac{\Ka(r)}{r}\,\P(V^2  < r)\, \Big
(1 +  O\,\Big(\frac{\ln{|\ln r|}}{|\ln r|}\Big) \Big). \eqno
(1.11)
$$
Moreover, for any $k\ge 2$
$$
\ep(r) = \sum_{m = 1}^k p_m(\bar s)\,|\ln{r}|^{-m}
+O\,((\ln{|\ln{r}|})^k/|\ln{r}|^{k+1}).\eqno (1.12)
$$ }

\medskip

Note that the remainders in (1.10) -- (1.12) are optimal.

\medskip

{\bf Corollary 1.} {\it For any $k\ge 2$ as $ r\to 0$
$$
\frac{8 \pi^2\, r}{\ln^2 r}\,|\log{\P(V^2 < r)}| = 1 + \sum_{m =
1}^k p_m(\bar s)\,|\ln{r}|^{-m} +
O\,((\ln{|\ln{r}|})^k/|\ln{r}|^{k+1}).
$$
In particular $($see $(1.8))$,
$$
\frac{8 \pi^2\, r}{\ln^2 r}\,|\log{\P(V^2 < r)}| = 1 + 2\,(\bar s
- 1)\,|\ln{r}|^{-1} + \bar s^2\,|\ln{r}|^{-2} +
O\,((\ln{|\ln{r}|})^2/|\ln{r}|^{3}).\eqno (1.13)
$$
}

\medskip

Note that (1.13) for $b = \d=-1/2$ refines [5,  Example 5.4 (d=2)]
and [14, Example 2] (simultaneously correcting the constant in
[14, (1.20)]).

\medskip

{\bf Corollary 2.} {\it Let $($see $(1.1))$
$$
\tilde V^2 = \pi^{-4}\sum\limits_{i,j\ge 1} \Big(((i+b)^2 -
q)\,((j + \d)^2 - \tau)\Big)^{-1} \,\xi_{i j}^2,
$$
where $q < 1 + b,\ \tau <  1 + \d$.

Then
$$
\P(\tilde V^2 < r)\sim C_{dist.}\,  \P(V^2 < r),\ r\to 0.
$$
Here
$$
C_{dist.} = \left(\frac{\Ga(1 + b + \sqrt{q})\,\Ga(1 + b
-\sqrt{q})\,\Ga(1 + \d + \sqrt{\tau})\,\Ga(1 + \d - \sqrt{\tau})}
{\Ga^2(1 + b)\,\Ga^2(1 + \d)}\right)^{1/2}.
$$
}

\medskip

In the conclusion of this section we present the results of
application of Theorem 1 to the classical cases mentioned earlier.
Here, the function $ \ep (r) = \ep (r, \bar s)  $ is defined in
(1.9), $ C_e $ is the Euler constant, and $\Ka(r) = (8 \pi^2\,
r)^{-1}\,\ln^2 r$.

\medskip



\medskip

{\bf Theorem 2.}{ Set  $\tilde s = \ln{|\ln{r}|} +  2\,C_e -
3\,\ln{\pi}$. Then as $r\to 0$

1)  The Brownian sheet $(\bar s = \tilde s + 3\,\ln{2})$
$$
\P(V^2\big|_{b = \d = -1/2} < r) =  \frac {e^{- \Ka(r)\,(1 +
\ep(r))}} {\sqrt {\pi\,\Ka(r)}}\,\Big (1 +  O\,\Big(\frac{\ln{|\ln
r|}}{|\ln r|}\Big) \Big).
$$

2) The Brownian pillow ($\bar s = \tilde s  -\ln{2}$)
$$
\P( V^2\big|_{b = \d = 0} < r) = r^{1/8}\,(2\,\Ka(r))^{-5/8}\,
\frac {e^{- \Ka(r)\,(1 + \ep(r))}} {\sqrt {\pi}}\,\Big (1 +
O\,\Big(\frac{\ln{|\ln r|}}{|\ln r|}\Big) \Big).
$$

3) The Kiefer field ($\bar s = \tilde s + \ln 2$)
$$
\P(V^2\big|_{b = -1/2,\,  \d = 0} < r) = 2^{-1/4}\, \frac {e^{-
\Ka(r)\,(1 + \ep(r))}} {\sqrt {\pi\,\Ka(r)}}\,\Big (1 +
O\,\Big(\frac{\ln{|\ln r|}}{|\ln r|}\Big) \Big).
$$
}


\bigskip

{\bf 2. Lemmas.}

\bigskip

 For $ b> -1 $ and $ \d> -1 $ denote
$$\gathered
\la(t) = (t + b)^{-1},\quad t\ge 1; \quad    \la_n = \la(n),\\
u(x) = \ln\Big(|\Ga(1 + \d + i\,x)|^2/\Ga^2(1 + \d)\Big),\\
J(\ga) = J(\ga;\,b, \d) = \sum\limits_{i\ge 1}
u(\la_i\,\ga).\endgathered\eqno (2.1)
$$

\medskip

{\bf Lemma 1.} {\it If\ $ \ga \to \infty $ then
$$\gathered
J(\ga) = -\pi\,\ga\,\ln\ga  + B_1\,\pi\,\ga  - B_2\,\ln\ga -
 B_3 + O\,(\ga^{-1}),\\
\ga\,J'(\ga) = -\pi\, \ga\,\ln\ga + (B_1 - 1)\,\pi\,\ga -  B_2 + O\,(\ga^{-1}),\\
\ga^2\,J''(\ga) =  -\pi\, \ga\ + B_2 +
O\,(\ga^{-1}),\endgathered\eqno (2.2)
$$
where $ B_i $ are defined in $ (1.2) $. }

\medskip

{\bf The proof of Lemma 1.}

\medskip

From  [1, (8.341.1)] provided $z = 1 + \d + i\,x$ it follows that
$$
u(x) = v_\d + (1/2 + \d)\,\log {(x^2 + (1 + \d)^2)} - 2\,x
\,(\arctan{x_\d} + 1/x_\d) + 2\,R(x),\eqno (2.3)
$$
with
$$
\gathered v_\d = \log{\frac{2\pi}{\Ga^2(1 + \d)}},\quad x_\d =
x/(1+\d),\\  R(x) = \int_0^\infty \Big(\frac{1}{2} - \frac{1}{t} +
\frac{1}{e^t-1}\Big)\,e^{-t\,(1+\d)}\,\cos{t
x}\,\frac{dt}{t}.\endgathered\eqno (2.4)
$$

Hence, $ u(x) = v(x) + g(x)$ and
$$
\gathered v(x) =  v_\d + (1 +2\,\d)\,\log x -\pi\,x,\\
g(x) = (1/2 + \d)\,\log{(1 + x^{-2}_\d)} - 2\,x \,(\arctan{x_\d} -
\pi/2 + 1/x_\d) + 2\,R(x).\endgathered\eqno (2.5)
$$

Note (see [1, 8.343.1]) that
$$
R^{(k)}(x)= O\,(x^{-k-2}),\quad g^{(k)}(x)= O\,(x^{-k-2})\quad
(k\ge 0),\quad x\to\infty.\eqno (2.6)
$$

Using the Euler -- Maclaurin  formula of the second order, we find
$$
J(\ga) = \int_1^\infty u(\ga\,\la(t))\,dt +
\frac{1}{2}\,u(\ga\,\la_1) + \Dl(\ga) + \bar\Dl(\ga),\eqno (2.7)
$$
where
$$
\Dl(\ga) = \sum\limits_{j\ge 1} \int\limits_0^1
\frac{t-t^2}{2}\,(v(\ga\,\la(t+j)))_{tt}''\,dt,\quad \bar\Dl(\ga)
= \sum\limits_{j\ge 1} \int\limits_0^1
\frac{t-t^2}{2}\,(g(\ga\,\la(t+j)))_{tt}''\,dt.
$$
In this case by (2.1), (2.5) и (2.6)
$$
|\bar\Dl(\ga)|\le \frac{1}{8}\, \int\limits_1^\infty
|(g(\ga\,\la(t)))_{tt}''|\,dt\le A_1\,\ga^{-1}\eqno (2.8)
$$
with the constant
$$
A_1 = \frac{1}{8}\, \int\limits_0^\infty (2\,y\,|g'(y)| +
y^2\,|g''(y)|) \,dy <\infty.
$$

Put
$$\gathered
s(t) = -(1 +2\,\d)\,\log t - \ga\,\pi\,t^{-1},\\
U_n = -\frac{1}{2}\,s(n + b) + \sum\limits_{j= 1}^n s(j + b) -
\int\limits_1^n s(t + b)\,dt.\endgathered\eqno (2.9)
$$

According to (2.5), $ v(\ga\,\la(t)))_{tt}'' = s''(t + b)$, and
therefore, by the Euler -- Maclaurin  formula  again,
$$
\Dl(\ga) = \sum\limits_{j\ge 1} \int\limits_0^1
\frac{t-t^2}{2}\,s''(t+j)\,dt = -\frac{1}{2}\,s(1 + b ) + \lim
U_n.\eqno (2.10)
$$
Straightforward calculations show that
$$
\Dl(\ga) = \ga\,\pi\,(A_2(b) - A_3(b)) + (1 +2\,\d)\,( A_4(b) -
A_5(b)),\eqno (2.11)
$$ 
where
$$\gathered
A_2(b) = (2\,(1 + b))^{-1} - \ln{(1 + b)}, \\ A_4(b) = 1 + b - (1/2 + b)\,\ln{(1 + b)},\\
A_3(b) = \lim\Big(\sum\limits_{j= 1}^n \frac{1}{j + b} - \ln
n\Big),\\
A_5(b) = \lim\Big(\sum\limits_{j=1}^n \ln{(j + b)} +n -(n + 1/2 +
b)\,\ln n\Big). \endgathered\eqno (2.12)
$$

The limits $ A_3 (b) $ and $ A_5 (b) $ may be calculated. So,
$$
A_3(b) = A_3(0) - b\,\sum\limits_{j\ge 1} 1/(j\,(j + b)) = -\Psi(1
+ b)\eqno (2.13)
$$
(see [1, 8.362.1,\, 8.366.1,\,8.366.2,\, 8.367.1,\,8.367.2]).

Similarly, taking into account [1, 8.322] (and the notation
(2.3)), we can show that
$$
A_5(b) = A_5(0) + b\, C_e + \sum\limits_{j\ge 1}\Big(\ln{(1 +
b/j)} - b/j\Big) = \frac{1}{2}\,\log{\frac{2\pi}{\Ga^2(1 + b)}} =
\frac{1}{2}\, v_b.\eqno (2.14)
$$

\medskip

Continue to analyze  the equality (2.7). From (2.1) and (2.3) --
(2.6) it follows that
$$
u(\ga\,\la_1) = v_\d + (1 +2\,\d)\,\log {(\ga/(1 + b))}
-\pi\,\ga/(1 + b) + O\,(\ga^{-2}), \ \ga\to\infty.\eqno (2.15)
$$
Further,
$$
\int_1^\infty u(\ga\,\la(t))\,dt = \ga\,\int\limits_0^{\ga\,\la_1}
u(x) \frac{dx}{x^2} =  -\frac{u(\ga\,\la_1)}{\la_1} +
\ga\,\int\limits_0^{\ga\,\la_1} u'(x)\,
\frac{dx}{x}.
\eqno (2.16)
$$

We have (see (2.3)),
$$
u'(x) =  - 2\,\arctan{x_\d} - \frac{1}{1 + \d}\,\frac{x_\d}{1 +
x^2_\d} + 2\, R'(x),
$$
and hence
$$\gathered
\int\limits_0^{\ga\,\la_1} u'(x)\, \frac{dx}{x} = -
2\,\int\limits_0^{\ga\,\la_1} \arctan{x_\d}\, \frac{dx_\d}{x_\d} -
\frac{1}{1 + \d}\,\arctan{(\ga\,\la_1/(1 + \d))} -\\
2\,\int\limits_{\ga\,\la_1}^\infty R'(x)\, \frac{dx}{x} +
2\,\int\limits_0^\infty R'(x)\, \frac{dx}{x}.
\endgathered\eqno (2.17)
$$
In this case (see (2.4) [1, \ 8.341.1] and (2.6)),
$$\gathered
\int\limits_0^\infty R'(x)\,\frac{dx}{x} = -\int_0^\infty
\Big(\frac{1}{2} - \frac{1}{t} + \frac{1}{e^t-1}\Big)\,e^{-t\,(1 +
\d)}\,\int_0^\infty\sin{t x}\,\frac{dx}{x}\ dt\\ =  \Big(\Psi(1 +
\d) - \log{(1 + \d)} + \frac{1}{2\,(1 + \d)}\Big)\,\frac{\pi}{2},\\
\int\limits_0^{\ga\,\la_1} \arctan{x_\d}\, \frac{dx_\d}{x_\d} =
\frac{\pi}{2}\,\log{\bar\ga} + \frac{1}{\bar\ga} + \,O(\ga^{-3}),\\
\arctan{(\ga\,\la_1/(1 + \d))} = \frac{\pi}{2} - \frac{1}{\bar\ga}
+ \,O(\ga^{-3}),
\endgathered\eqno (2.18)
$$
with $\bar\ga = \ga\,\la_1/(1 + \d) = \ga/((1 + \d)(1 + b))$.

\medskip

Combining (2.7) (2.8) (2.11) -- (2.18), we get the first relation
in (2.2).

\medskip

To obtain the other assertions in (2.2) one can differentiate the
right-hand  side of (2.7) in $\ga$ and use (2.6), (2.11) (2.16).

\medskip

Lemma 1 is proved.

\medskip

Further, let us give one more result, which will used in what
follows.

Let $ d> 0, \ 0 <\ep <d / e $. Consider the equation

$$
\ln{(d\, y)}/y = \ep,\ y>e/d, \eqno (2.19)
$$

and examine its solution $ y = y (\ep) $ for small positive $\ep$.

\medskip

{\bf Lemma 2.} {\it The solution of the equation $ (2.19) $ for
all $\ep $ small enough can be written as absolutely convergent
series
$$
y = \frac{\xi}{\ep}\, \left(1 + \sum_{j, l\ge 0}
a_{jl}\,s_d^{j+1}\,\xi^{-l-j-1}\right),\quad \xi = -\ln{\ep},\ s_d
= \ln{d\,\xi},
\eqno (2.20)
$$
where the coefficients $ a_ {jl} $ are defined in $(1.6)$.}

\medskip

Presumably, the representation (2.20) may be deduced from Theorem
2 in [3]. We provide a shorter proof below.

\medskip 

{\bf Proof of Lemma 2.}

\medskip

Put (see (2.20))
$$
y = \frac{\xi}{\ep}\,(1+v), \quad g(v) =  \ln{(1+v)}/v, \quad \si
= 1/\xi,\ \tau = s_d/\xi.\eqno (2.21)
$$

Then (2.20) is equivalent to the equation
$$
v - \si\, \ln{(1+v)} = v\,(1-\si g(v)) = \tau. 
$$
From this and from Lagrange's formula it follows that
$$
v  = \sum_{j\ge 1} d_j\,\tau^j, \eqno (2.22)
$$
for all $\ep$ small enough with
$$
d_j = \frac{1}{j!}\,\lim_{t\to 0}
\frac{d^{j-1}}{dt^{j-1}}\left(1-\si g(t)\right)^{-j}.
$$

We have,
$$
\left(1-\si g(t)\right)^{-j} = 1 + \sum_{l\ge 1} C^l_{j+l-1}
\si^l\, g^l(t),
$$
and hence
$$
d_{j+1} = \frac{1}{(j+1)!}\,\sum_{l\ge 0}\si^l\,
C^l_{j+l}\,\lim_{t\to 0}{ (g^l(t))^{(j)}},\ j\ge 0.\eqno (2.23)
$$

Taking into account that according to the definition of the
Stirling numbers of the first kind,
$$
\lim_{t\to 0}{ (g^l(t))^{(j)}} =  s(j+l, l)/C^l_{j+1},
$$
we obtain the first assertion of Remark 1. The second statement,
i.e. the equality (1.6),  follows from (2.21) -- (2.23) and [12,
Ch. VI, Lemma 3]. This also implies
$$
|a_{jl}|\le \frac{1}{j+1} \,C^l_{j+l}\,\frac{1}{2}
\,C^j_{j+l-1},
$$
and thus the estimate (1.9) holds.

Lemma 2 and, simultaneously, remark 1 are completely proved.

\medskip

Equation (2.20), obviously , can be rewritten in the form
$$
y = \frac{\xi}{\ep}\, \left(1 +\sum_{m\ge 1} \pi_m(s_d)\,
\xi^{-m}\right), \eqno (2.24)
$$
where the polynomials $ \pi_m (s) $ are defined in (1.5), and also
can be calculated using the recursive relation (1.4).

Check out the latter fact. Set $ f = f(z) = \sum_{m\ge 1}
\pi_{m+1}(s)\, z^m.$

From (2.19) and (2.24) the formal equality $f = \ln{(1 + s\,z +
z\,f)}$ follows. Differentiating it in $ z $, we find $(1 + s\,z
+z\,f)\,f' = s + (z\,f)'.$ Equating the coefficients of the
identical powers at $ z $ from the lefthand side and  from the
righthand side, we get
$$
\pi_{k+1} = \pi_k - \frac{1}{k}\, \sum_{i= 1}^{k-1}
i\,\pi_{k-i}\,\pi_{i+1} ,\ k\ge 1.
$$
Taking into account the fact that the sum is equal to $ \sum_{i=
0}^{k-2} (k-i-1)\,\pi_{k-i}\,\pi_{i+1}$, we obtain (1.4).

\medskip

{\bf Lemma 3.} {\it Let $ \xi_j $ and $ \la_j, \, j \ge 1 $ denote
independent standard Gaussian random variables and positive
numbers, respectively.

Assuming that $\sum\limits_{j\ge 1} \la_j < \infty$, put $S =
\sum\limits_{j\ge 1} \la_j\,\xi_j^2$ and for $h
> 0$ denote
$$
L(h) = -\frac{1}{2}\,\sum\limits_{j\ge 1} \ln{(1+2 h \la_j)},\
\tau^2(h)  = h^2\,L''(h). 
$$
Then for any $0<r<\E S$
$$
\P\big(S<r\big) =  \frac {1} {\tau(h)\sqrt {2\pi}}\, e ^ {L(h) +
h\,r} \left (1 + \theta_1\,\tau ^ {-2}(h) \right),
$$
$$
\frac{d}{dr} \P\big(S<r\big) =h\,\P\big(S<r\big)\,\left (1 +
\theta_2\,\tau ^ {-2}(h) \right),
$$
where $h$ is the unique solution of the equation
$$
L '(h) + r = 0, \eqno (2.25)
$$
and  $|\theta_i|$ is bounded by a constant}.

\medskip

Recall that
$$
L(h) + h\,r = \inf_{u>0}{(L(u) + u\,r)}. \eqno (2.26)
$$

\medskip

Lemma 3 follows from [13, Corollary 1].

\bigskip

{\bf 3. Proofs.}

\medskip

It is known (see [1, 8.326.1]) that
$$
\frac{\Ga^2(1 + \d)}{|\Ga(1 + \d + i\,x)|^2} = \prod_{j\ge
1}\Big(1+ \frac{x^2}{(j + \d)^2}\Big).
$$

Denoting for convenience $ h $ by $ \ga ^ 2/2 $, we get from here
(see the notation in (1.1) and (2.1))
$$
L(h)= \log{\E\,e^{-h V^2}} = \frac{1}{2}\, J(\pi^{-2}\,\ga; b, \d)
= \frac{1}{2}\,I(\ga).\eqno (3.1)
$$

\medskip

We set (see (1.3))
$$
\bar I(\ga) = - \pi^{-1}\,\ga\,\ln\ga +  B\,\ga -   C\,\ln\ga - D.
\eqno (3.2)
$$

\medskip

By Lemma 1, as $\ga\to\infty$,
$$
\gathered
I(\ga) = \bar I(\ga)  + O\,(\ga^{-1}),\\
$$
\ga\,I'(\ga) = \ga\,\bar I'(\ga)  + O\,(\ga^{-1}) =
-\pi^{-1}\,\ga\,\ln\ga + (B - \pi^{-1})\,\ga -  C  +
O\,(\ga^{-1}),\\
\ga^2\,I''(\ga) = \ga^2\,\bar I''(\ga)  + O\,(\ga^{-1})  =
-\pi^{-1}\,\ga +  C  + O\,(\ga^{-1}).\endgathered\eqno (3.3)
$$

\medskip

Now consider the equation (2.25) or (see (3.1))
$$
\frac{1}{2}\,I'(\ga) + \ga\,r = 0.\eqno (3.4)
$$

From (3.3) it follows, in particular, that the solution of the
equation (3.4) is such that
$$
\ga\sim \ga_0 = \ga_0(r) = \frac{|\ln r|}{2 \pi r},\ r\to 0.\eqno (3.5)
$$

We have (see (3.1), (3.2) and (2.26)) for positive $\ep $ small
enough
$$\gathered
2\,\inf_{u>0}{(L(u) + u\,r)} = \min_{1-\ep<s/\ga_0<1+\ep}{(I(s) +
s^2\,r)} \\=
\min_{1-\ep<s/\ga_0<1+\ep}{(\bar I(s) + s^2\,r)} + O\,(\ga_0^{-1})\\
= \bar I(\bar \ga) + \bar \ga^2\,r + O\,(\ga_0^{-1}),\ r\to 0,
\endgathered\eqno (3.6)
$$
where $\bar \ga = \bar \ga(r)\to\infty$ as $r\to 0$ is the
solution of the equation
$$
\frac{1}{2}\,\bar I'(\bar \ga) + \bar \ga\,r =0.\eqno (3.7)
$$

Put (see (1.3) and (3.2))
$$
\bar I_0(\ga) = \bar I(\ga) + C\,\ln \ga,\eqno
(3.8)
$$
and let  $\tilde \ga$ ($\tilde \ga\to\infty,\ r\to 0$) satisfy the
condition
$$
\frac{1}{2}\,\bar I_0'(\tilde \ga) + \tilde \ga\,r =0,\quad \text{
or }\quad \frac{\log\tilde \ga + 1 - \pi\,B}{2\,\pi\,\tilde \ga} =
r. \eqno (3.9)
$$

Now we compare the solutions of equations (3.7) and (3.9). Note
preliminary (see (3.5)) that $\bar \ga\sim\tilde \ga\sim \ga_0$ as
$ r \to 0 $.

We have by (3.2) and (3.3)),
$$\gathered
0 = \tilde \ga\,\bar I'(\bar \ga) - \bar \ga\,\bar I_0'(\tilde
\ga) = \tilde \ga\,(\bar I'(\bar \ga) - \bar I_0'(\bar \ga)) +
\tilde \ga\,(\bar I_0'(\bar \ga) - \bar I_0'(\tilde \ga)) - (\bar
\ga - \tilde \ga)\,\bar I_0'(\tilde \ga)\\
= -\tilde \ga\,\frac{ C}{\bar \ga} + \frac{\bar \ga - \tilde
\ga}{\tilde \ga}\,\Big(\tilde \ga^2\,\bar
I_0''(\hat\ga)\Big|_{\hat\ga\in {(\bar \ga, \tilde \ga})} - \tilde
\ga\,\bar I_0'(\tilde \ga)\Big)\\ = (-1+o\,(1))\, C + ( \pi^{-1} +
o\,(1))\, \frac{\bar \ga - \tilde \ga}{\tilde \ga}\,\chi(r),
\endgathered$$
where $ \chi(r) = \ga_0\,\ln \ga_0$.

This implies
$$
\frac{\bar \ga - \tilde \ga}{\tilde \ga} = O\,( |C|/\chi(r)),\
r\to 0. \eqno (3.10)
$$

Next, by (3.9)
$$
\gathered \bar I(\bar\ga) + \bar\ga^2\, r = - C\,\ln \bar\ga +
(\bar I_0(\bar\ga) + \bar\ga^2\, r)\\ = - C\,\ln \bar\ga + (\bar
I_0(\tilde\ga) + \tilde\ga^2\, r)+ \frac{1}{2}\,(\bar \ga - \tilde
\ga)^2\,\Big(2\,r +\bar I_0''(\hat\ga)\Big|_{\hat\ga\in {(\bar
\ga, \tilde \ga)}}\Big). 
\endgathered
$$

Therefore, taking into account (3.9) and (3.10), we obtain for
$r\to 0$
$$
\bar I(\bar\ga) + \bar\ga^2\, r = - C\,\ln \tilde\ga +  (\bar
I_0(\tilde\ga) + \tilde\ga^2\, r) + O\,(1/\chi(r)),
\eqno (3.11)
$$
and similarly
$$
\bar\ga^2\,\bar I''(\bar\ga) - \bar\ga\,\bar I'(\bar\ga) = 2\,C +
(\bar\ga^2\,\bar I_0''(\bar\ga) - \bar\ga\,\bar I_0'(\bar\ga)) =
\tilde\ga^2\,\bar I_0''(\tilde\ga) - \tilde\ga\,\bar
I_0'(\tilde\ga) + O\,(1). \eqno (3.12)
$$

\medskip

Now let us examine the behavior of $\tau^2(h) = h^2\,L''(h)$ as
$h\to\infty$.

It is obvious (see (3.1) and (3.3)) that
$$
\tau^2(h)\Big|_{h = \ga^2/2} =  \frac{1}{8}\,(\ga^2\,I''(\ga) -
\ga\,I'(\ga)) = \frac{1}{8}\,(\ga^2\,\bar I''(\ga) - \ga\,\bar
I'(\ga)) +  O\,(\ga^{-1}), \ \ga\to\infty.
$$
Therefore, if $\ga$ and $\tilde \ga$ satisfy the equations (3.7)
and (3.9), respectively (and $ h = \ga ^ 2/2 $), then
$$
\tau^2(h)\asymp \chi(r),\ r\to 0, \eqno (3.13)
$$
and, taking into account (3.10), (3.12),
$$ \tau^2(h) = \frac{1}{8}\,(\tilde\ga^2\,\bar I_0''(\tilde\ga) -
\tilde\ga\,\bar I_0'(\tilde\ga)) +  O\,(1),\ r\to 0.\eqno
(3.14)
$$

Denote
$$
Q(r) = -\frac{1}{2}\,(\bar I_0(\tilde\ga) + \tilde\ga^2\, r),\quad
q(r) = \frac{1}{8}\,(\tilde\ga^2\,\bar I_0''(\tilde\ga) -
\tilde\ga\,\bar I_0'(\tilde\ga)).
$$

By (3.2), (3.8) и (3.9)
$$
Q(r)  = \frac{1}{2}\,(\tilde\ga^2 \,r\ -  \pi^{-1}\,\tilde \ga +
D),\quad q(r) = \frac{1}{8}\,(2\,\tilde \ga^2\,r\ -
\pi^{-1}\,\tilde \ga),
\eqno (3.15)
$$
and, hence (see (3.6) and (3.11))
$$
\inf_{u>0}{(L(u) + u\,r)} = - Q(r) - \frac{C}{2}\,\ln \tilde\ga +
O\,(1/\ga_0),\ r\to 0.
$$

\medskip

This relation and Lemma 3 with $ S = V ^ 2 $, as well as relations
(2.26), (3.13) -- (3.15), imply
$$\gathered
\P\big(V^2 < r\big) = \big(2\pi\,\tilde\ga^{
C}\,q(r)\big)^{-1/2}\,
e ^ {-Q(r)} \left (1 + O\,(\ga^{-1}_0) \right),\\
\frac{d}{dr} \P\big(V^2 <r\big) = \frac{\tilde
\ga^2}{2}\,\P\big(V^2 < r\big)\, \left (1 + O\,(\ga^{-1}_0)
\right),\ r\to 0,
\endgathered \eqno  (3.16)
$$
provided  $\tilde\ga$ satisfies the equation (3.9) and $\ga_0$ is
defined in (3.5).

\medskip

Continue our reasoning.

Let $y = 2\,\pi\,\tilde\ga$, $d =  \pi^{-1}\,e^{1 - \pi\,B}/2$,\
$r =\ep$. Then (see (3.5) and (3.8)), the equation (3.9) has the
form (2.19).

Hence, using Lemma 2 ((2.24)), we obtain
$$\gathered
2\tilde \ga\,r =  \pi^{-1}\, \xi\,  \Big(1 + \sum\limits_{m\ge 1}
\bar\pi_m\, \xi^{-m}\Big), \quad \xi =
|\ln{r}|,\\
\bar \pi_m = \pi_m(\bar s),\, \quad \bar s = \ln \xi + 1 - \pi\, B
- \ln{( 2\, \pi)}.
\endgathered \eqno  (3.17)
$$

Put $\tilde\pi_m = \sum_{l=0}^m \bar \pi_l\,\bar \pi_{m-l}\quad
(\bar \pi_{0} = 1)$. We have,
$$
\tilde\pi_1 = 2\,\bar s, \quad \tilde\pi_2 = \tilde\pi_3 = 2\,\bar
s + \bar s^2, \eqno (3.18)
$$
and (see (1.4))
$$
\tilde\pi_m = -\frac{2}{m-2}\,(\bar \pi_m - (m-1)\,\bar
\pi_{m-1}),\quad m\ge 3.\eqno  (3.19)
$$

By (3.17)
$$
(2\tilde \ga\,r)^2 = \pi^{-2}\, \xi^2\,  \Big(1 +
\tilde\ep(r)\Big),\quad \tilde\ep(r) = \sum\limits_{m\ge 1}
\tilde\pi_m\, \xi^{-m}.\eqno (3.20)
$$

Using (3.15), (3.17) and (3.20), we find
$$\gathered
Q(r) = \frac{1}{2}\, D + \frac{ \pi^{-2}\, \xi^2}{8 r}\,\Big(1
+\ep(r)\Big),\quad \ep(r) = \sum\limits_{m\ge 1}
(\tilde\pi_m - 2\, \bar \pi_{m-1})\, \xi^{-m},\\
q(r) = \frac{\xi^2}{16\, \pi^{2}\,  r}\,\Big(1
+\ep_q(r)\Big),\quad \ep_q(r) = \sum\limits_{m\ge 1} (\tilde\pi_m
-  \bar \pi_{m-1})\, \xi^{-m}.
\endgathered
\eqno (3.21)
$$

The relations (3.16), (3.20) and (3.21) lead to the following
result.

\medskip

{\bf Proposition 1.} {\it Let
$$
\Ka(r) = \frac{\ln^2 r}{8\, \pi^2\, r},\quad K = \Big(\pi\,2^{
C/2}\,e^{ D}\Big)^{-1/2}.
$$
Then as $r\to 0$
$$\gathered
\P(V^2 < r) = K\,r^{ C/4}\,\Ka^{-( C + 2)/4}(r)\, e ^ {-
\Ka(r)\,(1+\ep(r))} \Big (1 + \d(r) + O\,(r/|\ln r|) \Big),\\
\frac{d}{dr} \P(V^2 < r) = K\,r^{ C/4 -1}\,\Ka^{(2 - C)/4}(r)\, e
^ {- \Ka(r)\,(1+\ep(r))}\, \Big (1 + \d(r) + \tilde\ep(r) +
O\,(r/|\ln r|) \Big),
\endgathered
$$
where
$$
\d(r) =  \Big((1 + \ep_q(r))\,(1 + \tilde\ep(r))^{
C/2}\Big)^{-1/2} -1.
$$
}

The function $\d(r) $ can be written as an absolutely convergent
series with the structure which is similar to the representation
of  $\tilde\ep(r)$ (see (3.20)). In particular,
$$
\d(r)  = \frac{1}{2}\,\Big(1  - (C + 2)\,\bar s\Big)/|\ln{r}|
 + O\,\Big(\frac{|\ln{|\ln{r}|}}{|\ln{r}|}\Big)^2.
\eqno  (3.22)
$$

The relations (1.10) and (1.11) of Theorem 1 follow from (3.22)
and Proposition 1 (recall that\ $ \Psi(1) = -C_e,\quad \Psi(1/2) =
- C_e - 2\,\ln 2$). The relation (1.12) is easily verified with
the help of (1.5)-- (1.7).

Corollary 2 follows from  the comparison theorem of [8] and the
well - known representation of gamma - function
$$
\frac{1}{\Ga(1 + z)} = e^{C_e\,z}\,\prod_{n\ge 1}(1 + z/n)\,
e^{-z/n}.
$$

\medskip

{\bf Acknowledgement.} The author is indebted to  professors M. A.
Lifshits and Ya.Yu. Nikitin (Saint-Petersburg State University)
for helpful discussions and valuable comments.

\bigskip

 \centerline {REFERENCES}

\medskip

1.  I.S. Gradshteyn, I.M. Ryzhik,  Tables of integrals, sums,
series and products, 5th ed. Moscow, Nauka, 1971 (in Russian).
English transl.: Table of integrals, series, and products. Corr.
and enl. ed. by Alan Je®rey. New York { London { Toronto: Academic
Press (Harcourt Brace Jovanovich, Publishers), 1980.

\medskip

2.  L. Comtet, Advanced Combinatorics (Reidel, 1974)).

\medskip

3.  L. Comtet, Inversion de $y^\a\,e^y$ et $y\,\log^\a{y}$ au
moyen des nombres de Stirling, Comptes Rendus de l'Academie
Scientifique, Ser. Mathematiques, 270(1970), 1085--8.

\medskip

4.  E. Cs\'aki, On small values of the square integral of a
multiparameter Wiener process,. In: Statistics and Probability,
Proc. of the 3rd Pannonian Symp. on Math. Stat. D.Reidel, Boston,
1982, 19–26.

\medskip

5.   J.A. Fill, F. Torcaso, Asymptotic analysis via Mellin
transforms for small deviations in L2-norm of integrated Brownian
sheets, Probab. Theory Relat. Fields 130 (2003), 259–288.

\medskip

6.   A. Karol’, A. Nazarov, Ya. Nikitin, Small ball probabilities
for Gaussian random fields and tensor products of compact
operators, Trans. AMS, 360 (2008), N3, 1443-–147.

\medskip

7.   A. Karol’,  A. Nazarov, Small ball probabilities for smooth
Gaussian fields and tensor products of compact operators,
Mathematische Nachrichten 287 (2014), 595–609,

\medskip

8.   W.V. Li, Comparison results for the lower tail of Gaussian
seminorms, J. Theor. Probab. 5(1) (1992), 1-–31.

\medskip

9.  A.I. Nazarov and Ya.Yu. Nikitin, Logarithmic L2-small ball
asymptotics for some fractional Gaussian processes, Theor. Probab.
and Appl., 49 (2004), 695–711.

\medskip

12.  V.V. Petrov, Sums of independent random variables}, Nauka,
Moscow 1972 (in Russian).

\medskip

13.  L. V. Rozovsky, On gaussian measure of balls in a hilbert
space, Theory Probab. Appl., 53(2) (2009), 357-364.

\medskip

14.  L. V. Rozovsky,  Small deviation probabilities for weighted
sums of independent random variables with a common distribution
that can decrease at zero fast enough, Statistics and Probability
Letters, 117 (2016), pp. 196--200.

\medskip

\end{document}